\documentclass[12pt]{amsart}
\numberwithin{equation}{section}
\newtheorem{theorem}{Theorem}
\numberwithin{theorem}{section}

\numberwithin{theorem}{section} \numberwithin{lemma}{section}

\numberwithin{definition}{section}
\newtheorem{corollary}{Corollary}
\numberwithin{corollary}{section}

\numberwithin{remark}{section}

\numberwithin{proposition}{section}
\def\b{\begin{equation}}
\def\e{\end{equation}}
\newcommand{\ignore}[1]{}
\normalsize
\date {February 2, 2006}
\thanks{}
\keywords{}
\begin{document}
\pagenumbering{arabic} \pagenumbering{arabic}\setcounter{page}{1}
\tracingpages 1
\title{improved  hardy and rellich inequalities on riemannian manifolds}
\author{Ismail Kombe and Murad \"Ozaydin}
\dedicatory {}
\address{Ismail Kombe, Mathematics Department\\ Dawson-Loeffler Science
\&Mathematics Bldg\\
Oklahoma City University \\
2501 N. Blackwelder, Oklahoma City, OK 73106-1493}
\email{ikombe@okcu.edu}
\address{Murad Ozaydin, Mathematics Department\\ University of
Oklahoma, Norman,  OK 73019-0315 } \email{mozyadin@math.ou.edu}

\begin{abstract}
In this paper we establish improved  Hardy and Rellich type
inequalities on Riemannian manifold $M$. Furthermore, we also
obtain sharp constant for the improved Hardy  inequality and
explicit constant for the Rellich inequality on hyperbolic space
$\mathbb{H}^n$.
\end{abstract}
\maketitle
\section{Introduction}
 The classical
Hardy inequality states that for $n\ge 3$
\begin{equation}
\int_{\mathbb{R}^n}|\nabla\phi(x)|^2dx\ge
\Big(\frac{n-2}{2}\Big)^2\int_{\mathbb{R}^n}
\frac{|\phi(x)|^2}{|x|^2}dx,
\end{equation}
where $\phi\in C_0^{\infty}( \mathbb{R}^n\setminus \{0\})$ and the
constant $(\frac{n-2}{2})^2$ is sharp. An extension of the Hardy's
inequality to  second-order derivative is the well known Rellich
inequality:
\begin{equation} \int_{\mathbb{R}^n}|\Delta \phi(x)|^2dx\ge
\frac{n^2(n-4)^2}{16}\int_{
\mathbb{R}^n}\frac{|\phi(x)|^2}{|x|^4}dx\end{equation}
 for all
$\phi\in C_c^{\infty}(\mathbb{R}^n\setminus \{0\}$) and $n\neq 2$,
where the constant $\frac{n^2(n-4)^2}{16}$   is sharp.
\medskip

 Over the last twenty years, there
has been a lot of research concerning Hardy and Rellich
inequalities on the Euclidean space $\mathbb{R}^n$ and, in
particular, sharp inequalities as well as their improved versions
which have attracted a lot of attention because of their
application to singular problems, e.g. \cite{Baras-Goldstein},
\cite {Brezis-Vazquez}, \cite{Garcia-Peral},
\cite{Vazquez-Zuazua}, \cite{Cabre-Martel},
\cite{Goldstein-Zhang}, \cite{Dupaigne}, \cite{Kombe1} and
references therein. In recent years, much attention has been paid
to Hardy and Rellich inequalities in Sub-Riemmannian spaces, e.g.
\cite{Kombe2} and references there in. In contrast, there is
considerably less literature for general Riemannian manifold. In
an interesting paper, Carron \cite{Carron} studied weighted
$L^2$-Hardy inequalities under
some geometric assumptions on the weight function $\rho$ and obtained, among other results, the following inequality:\\

\begin{equation}\int_M \rho^{\alpha}|\nabla \phi|^2dx\ge
\Big(\frac{C+\alpha-1}{2}\Big)^2\int_
M\rho^{\alpha}\frac{\phi^2}{\rho^2}dx\end{equation} where
$\alpha\in \mathbb{R}$, $C+\alpha-1>0$, $\phi\in
C_0^{\infty}(M-\rho^{-1}\{0\})$ and the weight function satisfies
$|\nabla\rho|=1$ and $\Delta\rho\ge \frac{C}{\rho}$ in the sense
of distribution.

Under these geometric assumptions, our first goal is to obtain
weighted Hardy  and Rellich type inequalities with remainder
terms.   We should  mention that Davies and Hinz \cite
{Davies-Hinz} studied $L^p$-Rellich type inequalities as well as
their higher order versions. In \cite {Grillo}, Grillo  obtained
Hardy, Rellich and Sobolev inequalities in the context of
homogeneous spaces. Recently, Barbatis \cite{Barbatis} obtained,
under a geometric assumption, an improved higher-order Rellich
inequality of the form
\begin{equation}
\int_{\Omega}\frac{|\Delta ^{m/2} u|^p}{d^{\gamma}}dx\ge
A(m,\gamma)\int_{\Omega}
\frac{|u|^p}{d^{\gamma+mp}}dx+B(m,\gamma)\int_{\Omega}V_i|u|^pdx,\end{equation}
where $u\in C_c^{\infty}(\Omega\setminus K)$, $d(x)=\text{dist}(x,
K)$, $\gamma\in \mathbb{N}$ and $V_i(x)$ involves some suitable
iterated logarithmic functions. Our improved Rellich inequality,
Theorem 2.3 below, has a different type of a remainder term
involving gradient $\nabla u$. Our method is different and
simpler. Let us recall that a complete Riemannian manifold $M$ is
said to be nonparabolic if there exists a symmetric positive
Green's function $G(x,y)$ for the Laplacian $\Delta$ on $L^2$
functions.  Recently, Li  and Wang \cite{Li-Wang} proved that,
among other results,  existence of a weighted Hardy type
inequality is equivalent to nonparabolicity.  Furthermore, they
obtained the following $L^2$-Hardy inequality

\[\int_M |\nabla\phi|^2dx\ge \int_M \frac{|\nabla
G(p,x)|^2}{4G^2(p,x)}\phi^2dx\] where $\phi\in C_c^{\infty}(M)$
and  $G(p,x)$ is the minimal positive Green's function defined on
$M$ with a pole at the point $p\in M$.
\medskip

 The second goal of this paper is to find
sharp versions of improved Hardy inequalities and an improved
Rellich inequality in the specific case of the hyperbolic spaces
$\mathbb{H}^n$. Both Hardy and Rellich inequalities in hyperbolic
spaces as well as the determination of sharp constants is new.

\section{ Improved Hardy and Rellich inequalities on Riemannian manifolds }
In the various integral inequalities below (Section 2 and Section
3), we allow the values of the integrals on the left-hand sides to
be $+\infty$. The following theorem is the first result of this
section.

\begin{theorem}
Let $M$ be a complete noncompact Riemannian manifold of dimension
$n>1$. Let $\rho$ be a nonnegative function on $M$ such that
$|\nabla\rho|=1$ and $\Delta\rho\ge \frac{C}{\rho}$ where $C>q-1$.
Then the following inequality hold
\[\Big (\int
_M|\phi|^p\Big)^{\frac{p-q}{p}}\Big(\int _M
|\nabla\phi|^p\Big)^{\frac{q}{p}}dx\ge
\Big(\frac{C+1-q}{p}\Big)^q\int _M \frac{|\phi|^p}{\rho^ q}dx\]
for all compactly supported smooth function $\phi\in
C_c^{\infty}(M\setminus \rho^{-1}\{0\})$, where  $1\le p<\infty$
and $0\le q \le p $. \proof The argument is a simple application
of the divergence theorem, as follows:\\ Let $Q(x)=\frac{\nabla
\rho}{\rho^{q-1}}$ then we have $\nabla\cdot Q\ge
\frac{(C+1-q)}{\rho^{q}}$. It is clear that
\[\nabla\cdot (|\phi|^p
Q(x))=p|\phi|^{p-1}\nabla|\phi|\cdot Q(x)+|\phi|^p\nabla\cdot
Q(x).\] Integrating the above formula, we get

\[(C+1-q)\int_M \frac{|\phi|^p}{\rho^q}dx\le p\int_M
\frac{|\phi|^{p-1}}{\rho^{q-1}}|\nabla \phi|dx.\] Applying
H\"older's inequality, we obtain the desired inequality
\begin{equation}
\Big (\int _M|\phi|^pdx\Big)^{\frac{p-q}{p}}\Big(\int _M
|\nabla\phi|^p\Big)^{\frac{q}{p}}dx\ge
\Big(\frac{C+1-q}{p}\Big)^q\int _M \frac{|\phi|^p}{\rho^ q}dx
\end{equation}
\qed
\end{theorem}
\medskip

We are now ready to prove an $L^p$-version of the inequality
(1.3).
\begin{theorem}
Let $M$ be a complete noncompact Riemannian manifold of dimension
$n>1$. Let $\rho$ be a nonnegative function on $M$ such that
$|\nabla\rho|=1$ and $\Delta\rho\ge \frac{C}{\rho}$ in the sense
of distribution where $ C>0$. Then the following inequality hold
\begin{equation}\Big (\int _M\rho ^{\alpha}|\nabla\phi|^p dx\ge
\Big(\frac{C+1+\alpha-p}{p}\Big)^p\int _M
\rho^{\alpha}\frac{|\phi|^p}{\rho^ p}dx\end{equation} for all
compactly supported smooth function $\phi\in
C_c^{\infty}(M\setminus \rho^{-1}\{0\})$, $ 1\le p<\infty$ and $
C+1+\alpha-p>0$.

\proof Let $\phi=\rho^{\gamma}\psi$  where $\psi\in
C_c^{\infty}(M)$ and $\gamma<0$. A direct calculation
 shows that
\[
|\nabla(\rho^{\gamma}\psi)|=|\gamma \rho^{\gamma-1}\psi\nabla
\rho+ \rho^{\gamma}\nabla \psi|.
\]
We now use the following inequality which is  valid for any $a, b
\in \mathbb{R}^n$ and  $p>2$,

\[
|a+b|^p-|a|^p\ge c(p)|b|^p+p|a|^{p-2}a\cdot b
\]
where $c(p)> 0$. This yields

\[\rho^{\alpha}|\nabla\phi|^p\ge |\gamma|^p \rho^{\gamma
p-p+\alpha}|\psi|^p+p|\gamma|^{p-2}\gamma \rho^{\alpha +\gamma
p+1-p}|\psi|^{p-2}\psi \nabla\rho\cdot \nabla\phi.\] Then
integration by parts gives

\[\begin{aligned}\int_M \rho^{\alpha}|\nabla\phi|^pdx &\ge |\gamma|^p \int_M \rho^{\gamma
p-p+\alpha}|\psi|^pdx-\frac{|\gamma|^{p-2}\gamma}{\alpha+\gamma
p-p+2} \int_M \Delta (\rho^{\alpha+\gamma p-p+2})|\psi|^p dx\\
&\ge (1-p)|\gamma|^p \int_M \rho^{\gamma p-p+\alpha}|\psi|^pdx
-\gamma|\gamma|^{p-2}(\alpha+C+1-p)\int_M \rho^{\gamma p-p
+\alpha}|\psi|^pdx. \end{aligned}\] We now choose
$\gamma=\frac{p-\alpha-C-1}{p}$ then we get the desired inequality
\[\int_M \rho^{\alpha}|\nabla\phi|^pdx \ge
\Big(\frac{C+1+\alpha-p}{p}\Big)^p\int_M \rho^{\alpha}
\frac{|\phi|^p}{\rho^p}dx.\] The theorem (2.2) also holds for
$1<p<2$ and in this case we use the following inequality
\[|a+b|^p-|a|^p\ge
c(p)\frac{|b|^2}{(|a|+|b|)^{2-p}}+p|a|^{p-2}a\cdot b
\]
where $c(p)> 0$ (see \cite {Lindqvist}). \qed
\end{theorem}

We now prove the following improved Hardy inequality which is
inspired by a recent work of Abdellaoui, Colorado and Peral \cite
{Peral}.

\begin{theorem} Let $M$ be n-dimensional complete noncompact Riemannian manifold and let
$\rho$ be nonnegative function such that $|\nabla\rho|=1$ and
$\Delta\rho\ge \frac{C}{\rho}$ in the sense of distribution where
$C>0$. Let $\Omega$ be a bounded domain with smooth boundary which
contains origin, $1<q<2$, $\alpha\in \mathbb{R}$, $C+\alpha-1>0$,
$\phi\in C_0^{\infty}(\Omega)$  then there exists a positive
constant $C_1=C_1(n, q, \Omega)$ such that the following
inequality is valid
\begin{equation}\int_{\Omega}\rho^{\alpha}|\nabla\phi|^2dx\ge
\Big(\frac{C+\alpha-1}{2}\Big)^2\int_{\Omega}\rho^{\alpha}\frac{\phi^2}{\rho^2}dx+C_1(\int_{\Omega}|\nabla\phi|^q
\rho^{q\alpha/2}dx)^{2/q}\end{equation} \proof
 Let $\psi\in C_c^{\infty}(M)$ then a
straight forward computation shows that

\[|\nabla\phi|^2-\nabla(\frac{\phi^2}{\psi})\cdot\nabla\psi=
\Big|\nabla\phi-\frac{\phi}{\psi}\nabla\psi\Big|^2.\]

Therefore
\[\begin{aligned}\int_{\Omega}\Big(|\nabla\phi|^2-\nabla(\frac{\phi^2}{\psi})\cdot\nabla\psi\Big)\rho^{\alpha}dx &=
\int_{\Omega}\Big|\nabla\phi-\frac{\phi}{\psi}\nabla{\psi}\Big|^2\rho^{\alpha}dx\\
&\ge
C_1\Big(\int_{\Omega}\Big|\nabla\phi-\frac{\phi}{\psi}\nabla\psi\Big|^q\rho^{q\alpha/2}
dx\Big)^{2/q}\end{aligned}\] where we used the Jensen's inequality
in the last step. Let us choose $\psi=\rho^{\beta}$ where
$\beta<0$. Then it is clear that

\[\begin{aligned}\int_{\Omega}\Big(|\nabla\phi|^2-\nabla(\frac{\phi^2}{\psi})\cdot\nabla\psi\Big)\rho^{\alpha}dx&=
\int_{\Omega}\rho^{\alpha}|\nabla\phi|^2dx +\frac{\beta}{\alpha+\beta} \int_{\Omega}(\frac{\Delta(\rho^{\alpha+\beta})}{\rho^{\beta}})\phi^2dx\\
&\le\int_{\Omega}\rho^{\alpha}|\nabla\phi|^2dx +
\beta(\alpha+\beta+C-1)\int_{\Omega}\rho^{\alpha-2}\phi^2dx.\end{aligned}\]
Therefore we have
\begin{equation}\int_{\Omega}\rho^{\alpha}|\nabla\phi|^2dx\ge
-\beta^2-\beta(\alpha+C-1)\int_{\Omega}\rho^{\alpha}\frac{\phi^2}{\rho^2}dx+
C_1\Big(\int_{\Omega}\Big|\nabla\phi-\frac{\phi}{\psi}\nabla\psi\Big|^q
\rho^{q\alpha/2}dx\Big)^{2/q}.\end{equation} Now we can use the
following elementary inequality : Let $1<q<2$ and $w_1, w_2\in
\mathbb{R}^n$ then the following inequality hold:

\begin{equation}c(q)|w_2|^q\ge |w_1+w_2|^q-|w_1|^q-q|w_1|^{q-2}\langle w_1,
w_2\rangle.\end{equation} Therefore by integration and using
successively the inequality (2.5), Young's and weighted
$L^p$-Hardy inequality (2.2),  we get

\begin{equation}\int_{\Omega}\Big|\nabla\phi-\frac{\phi}{\psi}\nabla\psi\Big|^q\rho^{q\alpha/2}
dx\ge C_1 \int _{\Omega}|\nabla\phi|^q\rho^{q\alpha/2}dx.
\end{equation}
Substituting (2.6)  into (2.4)  then we get
\[\int_{\Omega}\rho^{\alpha}|\nabla\phi|^2dx\ge -\beta^2-\beta(\alpha+C-1)\int_{\Omega}\rho^{\alpha}\frac{\phi^2}{\rho^2}dx+
C_1\Big(\int_{\Omega}|\nabla_{\mathbb{G}}\phi|^q
\rho^{q\alpha/2}dx\Big)^{2/q}.\] Now choosing
$\beta=\frac{1-\alpha-C}{2}$ then we obtain the desired inequality
\[\int_{\Omega}\rho^{\alpha}|\nabla\phi|^2dx\ge \Big(\frac{C+\alpha-1}{2}\Big)^2\int_{\Omega}\rho^{\alpha}\frac{\phi^2}{\rho^2}dx+
C_1\Big(\int_{\Omega}|\nabla\phi|^q
\rho^{q\alpha/2}dx\Big)^{2/q}.\] \qed
\end{theorem}

We now prove the following improved  Rellich inequality. In the
Euclidean case, our result improve a result of Davies and Hinz
\cite {Davies-Hinz}.

\begin{theorem}(Improved Rellich Inequality)
Let $M$ be a complete noncompact Riemannian manifold of dimension
$n>1$. Let $\rho$ be a nonnegative function on $M$ such that
$|\nabla\rho|=1$ and $\Delta\rho\ge \frac{C}{\rho}$ in the sense
of distribution where $C>0$. Then the following inequality hold
\[\begin{aligned}\int_{M}\rho^{\alpha}|\Delta \phi|^2dx &\ge
\frac{(C+\alpha-3)^2(C-\alpha+1)^2}{16}
\int_{M}\rho^{\alpha}\frac{\phi^2}{\rho^4}dx\\
&+\frac{(C+\alpha-3)(C-\alpha+1)}{2}C_1\int_M
|\nabla\phi|^q\rho^{\frac{q(\alpha-2)}{2}}dx\big)^{2/q}\end{aligned}\]
for all compactly supported smooth function $\phi\in
C_c^{\infty}(M\setminus \rho^{-1}\{0\})$, $\alpha<2$ and
$C+\alpha-3>0$. \proof A straight forward computation shows that
\begin{equation}\Delta \rho^{\alpha-2}\le (\alpha-2)(C+\alpha-3)\rho^{\alpha-4}.\end{equation}
Multiplying both sides of (2.7) by $\phi^2$ and integrating over
$M$,  we obtain
\begin{equation}\begin{aligned}(C+\alpha-3)(\alpha-2)\int_M \rho^{\alpha-4}\phi^2dx &\ge \int_M \rho^{\alpha-2}
\Delta(\phi^2)dx\\& = \int_M
\rho^{\alpha-2}(2|\nabla\phi|^2+2\phi\Delta\phi)dx.
\end{aligned}
\end{equation}
Therefore
\begin{equation}-2\int_M(\phi\Delta\phi)\rho^{\alpha-2}\ge 2\int_M
\rho^{\alpha-2}|\nabla\phi|^2dx-(\alpha-2)(C+\alpha-3)\int_M
\rho^{\alpha-4}\phi^2dx.\end{equation} After we apply weighted
Hardy and Cauchy-Schwarz inequalities, we obtain the following
plain weighted Rellich inequality
\[\int_{M}\rho^{\alpha}|\Delta \phi|^2dx \ge
\frac{(C+\alpha-3)^2(C-\alpha+1)^2}{16}
\int_{M}\rho^{\alpha}\frac{\phi^2}{\rho^4}dx.\]

Furthermore, let us apply Young's inequality to expression
$-2\int_M(\phi\Delta\phi)\rho^{\alpha-2}$ in (2.9) and we obtain
\begin{equation}-\int _{M}\rho^{\alpha-2}\phi\Delta\phi dx\le \epsilon\int
_{M}\rho^{\alpha-4}\phi^2dx+\frac{1}{4\epsilon} \int
_{M}\rho^{\alpha} |\Delta\phi|^2dx\end{equation} where
$\epsilon>0$. Substituting (2.10) into (2.9) and using the
improved Hardy inequality (2.3), we get
\[\frac{1}{4\epsilon}\int _M\rho^{\alpha}
|\Delta\phi|^2dx\ge
\Big(\frac{(C+\alpha-3)(C-\alpha+1)}{4}-\epsilon\Big)\int
_M\rho^{\alpha-4}\phi^2dx+C_1\big(\int_M
|\nabla\phi|^q\rho^{\frac{q(\alpha-2)}{2}}dx\big)^{2/q}.\] Since
$C+\alpha-3>0$ and $C-\alpha+1>0$ then we choose
$\epsilon=\frac{(C+\alpha-3)(C-\alpha+1)}{8}$.  Therefore  we
obtain the following improved Rellich inequality
\begin{equation}\begin{aligned}
\int _M\rho^{\alpha} |\Delta\phi|^2dx &\ge
\frac{(C+\alpha-3)^2(C-\alpha+1)^2}{16} \int_M\rho^{\alpha}\frac{\phi^2}{\rho^4}dx\\
&+\frac{(C+\alpha-3)(C-\alpha+1)}{2}C_1\int_M
|\nabla\phi|^q\rho^{\frac{q(\alpha-2)}{2}}dx\big)^{2/q}
.\end{aligned}\end{equation} \qed
\end{theorem}

\noindent\textbf{Uncertainty Principle Inequality.} The classical
uncertainty principle was developed in the context of quantum
mechanics by Heisenberg \cite {Heisenberg}. It says that the
position and momentum of a particle cannot be determined exactly
at the same time but only with an ``uncertainty". The harmonic
analysis version of uncertainty principle states that a function
on the real line and its Fourier transform can not be
simultaneously well localized. It has been widely studied in
quantum mechanics and signal analysis. There are various forms of
the uncertainty principle. For an overview we refer to Folland's
and Sitaram's paper \cite {Folland-Sitaram}.

The uncertainty principle on the Euclidean space $ \mathbb{R}^n$
can be  stated in the following way:
\begin{equation}
\Big(\int_{\mathbb{R}^n} |x|^2
|f(x)|^2dx\Big)\Big(\int_{\mathbb{R}^n} |\nabla f(x)|^2 dx\Big)\ge
\frac{n^2}{4} \Big(\int_{\mathbb{R}^n} |f(x)|^2 dx \Big)^2
\end{equation}
for all $f\in L^2( \mathbb{R}^n)$.

Using the Hardy type inequalities, we obtain the following
uncertainty principle type inequalities on Riemannian manifold
$M$.

\begin{corollary}($L^p$-Uncertainty type inequality) Let $M$ be a complete noncompact Riemannian manifold of dimension $n>1$. Let
$\rho$ be a nonnegative function on $M$ such that $|\nabla\rho|=1$
and $\Delta\rho\ge \frac{C}{\rho}$ in the sense of distribution
where $C>0$. Then the following inequality hold
\[\Big(\int_M \rho^q
\phi^qdx\Big)^{1/q}\Big(\int_M |\nabla\phi|^pdx\Big)^{1/p}\ge
\frac{C+1-p}{p}\int_M \phi^2dx\] for all compactly supported
smooth function $\phi\in C_c^{\infty}(M\setminus \rho^{-1}\{0\})$,
$ 1< p<\infty$, $\frac{1}{p}+\frac{1}{q}=1$ and $C+1-p>0$.
\end{corollary}
\medskip

\begin{corollary}(Improved $L^2$-Uncertainty type inequality) Let $M$ be a
complete noncompact Riemannian manifold of dimension $n>1$. Let
$\rho$ be a nonnegative function on $M$ such that $|\nabla\rho|=1$
and $\Delta\rho\ge \frac{C}{\rho}$ in the sense of distribution
where $C>0$. Then the following inequality hold
\[\Big(\int_M \rho^{\alpha}
\phi^2dx\Big)\Big(\int_M
\rho^{\alpha}|\nabla\phi|^2dx-C_1\Big(\int_M
|\nabla\phi|^2\rho^{\frac{\alpha q}{2}}\Big)^{2/q}\Big)\ge
\Big(\frac{C+\alpha-1}{2}\Big)^2\Big(\int_M \phi^2dx\Big)^2\] for
all compactly supported smooth function $\phi\in C_c^{\infty}(M)$,
$ 1<q<2$, $C+\alpha-1>0$ and $C_1>0$.
\end{corollary}
\medskip

\medskip
\section {Sharp improved Hardy and Rellich inequalities on hyperbolic space $\mathbb{H}^n $}

We will be using the Poincare conformal disc model for the
hyperbolic space $\mathbb{H}^n$. So the underlying space is
\[\mathbb{B}^n=\{x=(x_1, \cdots, x_n) \in \mathbb{R}^n |\, |x|<1\}\] in
$\mathbb{R}^n$equipped with the Riemannian metric obtained by
scaling the Euclidean metric with a factor of $
p:=\frac{2|dx|}{1-|x|^2}.$ Hence $\{pdx_i\}_{i=1}^{n}$ give an
orthonormal basis of the tangent space at $x=(x_1,\cdots, x_n)$ in
$\mathbb{B}^n$. The corresponding dual basis is
$\{\frac{1}{p}\frac{\partial}{\partial x_i}\}_{i=1}^{n}$, thus the
hyperbolic gradient is
\[\nabla_{\mathbb{H}^n}u=\frac{\nabla
u}{p}\] where $u\in C^1(\mathbb{B}^n)$ and $\nabla u$ is the usual
gradient. $\mathbb{H}^n$ is a contractible complete Riemannian
manifold with all sectional curvatures equal $-1$.  Geodesic lines
passing through the origin are the diameters of $ \mathbb{B}^n$
along with open arcs of circles in $\mathbb{B}^n$ perpendicular to
the boundary at $\infty$,  $\partial \mathbb{B}^n=S^{n-1}= \{x\in
\mathbb{R}^n: |x|=1\}$. It follows that the distance from $x\in
\mathbb{B}^n$ to the origin is \[ d=d_{\mathbb{H}}(0,x)=\log
(\frac{1+|x|}{1-|x|})\]; the hyperbolic volume element is:
\[dV=p^n(x) dx\] (where $dx$ is the usual Euclidean volume
element) and the  Laplace-Beltrami operator is given by
\[\Delta_{\mathbb{H}^n}u=p^{-n}\text{div}(p^{n-2}\nabla u)
\]
where $\nabla$ and $\text{div}$ denote the Euclidean  gradient and
divergence in $\mathbb{R}^n$, respectively.

Note that we have the following two relations for the distance
function $d= \log (\frac{1+|x|}{1-|x|})$ :

\[\begin{aligned}
 |\nabla_{\mathbb{H}^n}d|\, &=1,\\
 \Delta_{\mathbb{H}^n} d\, &\ge
\frac{n-1}{d}, \quad x \neq 0.\end{aligned}\]

Let us remark that the Poincare inequality with the Muckenhoupt
weight play an important role in the following theorem. We recall
that a weight $w(x)$ satisfies Muckenhoupt $A_p$ condition for
$1<p<\infty$ if there is a constant $C$ such that
\[\Big(\frac{1}{|B|}\int_{B}w(x)dx\Big)^{1/p}\Big(\frac{1}{|B|}\int_{B}w(x)^{-p'/p}dx\Big)^{\frac{1}{p'}}\le
C\] for all  balls $B$.  If $w(x)\in A_p$ then we have  $
w(x)^{-p'/p}\in A_{p'}$ where $p'$ is the dual exponent to $p$
given by $\frac{1}{p}+\frac{1}{p'}=1$.

Before we state  our first theorem, let us state a well known
result of Brezis and
 V\'azquez \cite {Brezis-Vazquez} in this direction. They proved that for a bounded domain $\Omega\subset
 \mathbb{R}^n$ there holds
 \begin{equation}
\int _{\Omega} |\nabla \phi(x)|^2dx\ge
\Big(\frac{n-2}{2}\Big)^2\int_{\Omega}
\frac{|\phi(x)|^2}{|x|^2}dx+\mu\big(\frac{\omega_n}{|\Omega|}\big)^{2/n}\int
_{\Omega} \phi^2dx,
 \end{equation}
where $\omega_n$  and $|\Omega|$ denote the $n$-dimensional
Lebesgue measure of the unit ball $B\subset \mathbb{R}^n$  and the
domain $\Omega$ respectively. Here $\mu$ is the first eigenvalue
of the Laplace operator in the two dimensional unit disk and it is
optimal when $\Omega$ is a ball centered at the origin. We now
prove a similar (weighted) improved Hardy inequality on hyperbolic
space $\mathbb{H}^n$.

\begin{theorem}
Let $\alpha\in \mathbb{R}$ and $\phi\in
C_0^{\infty}(\mathbb{H}^n\setminus \{0\})$. Then we have :
\[\int_{
\mathbb{B}^n }d^{\alpha}|\nabla \phi|^2p^{n-2}dx\ge
\Big(\frac{n+\alpha-2}{2}\Big)^2\int_ {\mathbb
{B}^n}d^{\alpha}\frac{\phi^2}{d^2} p^n dx+
c2^{n-2}\int_{\mathbb{B}^n }d^{\alpha}\phi^2 dx\] where $d=\log
(\frac{1+|x|}{1-|x|})$ is the distance from $x\in \mathbb{B}^n $
to the origin and $c>0$. Moreover, the constant
$(\frac{n+\alpha-2}{2})^2$ is sharp provided $n+\alpha-2>0$.
\end{theorem}
\proof Let $\phi=d^{\beta}\psi$ where $\beta\in
\mathbb{R}\setminus \{0\}$ and $\psi \in
C_0^{\infty}(\mathbb{H}^n\setminus \{0\})$. A direct calculation
shows that
\begin{equation}
d^{\alpha}|\nabla\phi|^2p^{n-2} =\beta^2d^{\alpha+2\beta-2}|\nabla
d|^2\psi^2p^{n-2}+ 2\beta d^{\alpha+2\beta-1}\psi p^{n-2}\nabla
d\cdot\nabla\psi+d^{\alpha+2\beta}|\nabla\psi|^2 p^{n-2}.
\end{equation}
It is easy to see that
\[|\nabla d|^2=p^2\]
and integrating (3.2) over $ \mathbb{B}^n$, we get
\begin{equation}
\begin{aligned}\int_{\mathbb{B}^n}d^{\alpha}|\nabla\phi|^2 p^{n-2}dx&=\int_{
\mathbb{B}^n}\beta^2 d^{\alpha+2\beta-2}\psi^2
p^ndx+\int_{\mathbb{B}^n} 2\beta d^{\alpha+2\beta-1}\psi
p^{n-2}\nabla d\cdot\nabla\psi dx\\&+\int_{\mathbb{B}^n}
d^{\alpha+2\beta}|\nabla\psi|^2 p^{n-2}dx.\\
\end{aligned}
\end{equation}
Applying integration by parts to the middle integral on the
right-hand side of (3.3), we obtain

\begin{equation}
\begin{aligned}
\int_{\mathbb{B}^n}d^{\alpha}|\nabla\phi|^2 p^{n-2}dx&=\int_{
\mathbb{B}^n}\beta^2 d^{\alpha+2\beta-2}\psi^2 p^ndx
 -\frac{\beta}{\alpha+2\beta} \int_{\mathbb{B}^n}\text{div}\big(p^{n-2}\nabla (d^{2\beta+\alpha})\big)dx \\
 &+\int_{\mathbb{B}^n}
d^{\alpha+2\beta}|\nabla\psi|^2 p^{n-2}dx.
\end{aligned}
\end{equation}
One can show that
\begin{equation}
\begin{aligned}
&-\frac{\beta}{\alpha+2\beta}
\int_{\mathbb{B}^n}\text{div}\big(p^{n-2}\nabla
(d^{2\beta+\alpha})\big)dx \\ =& -\beta (2\beta
+\alpha-1)\int_{\mathbb{B}^n}d^{2\beta+\alpha-2}p^n\psi^2dx-\beta
\int_{\mathbb{B}^n}d^{2\beta+\alpha-1}p^{n-2}\psi^2 (\Delta d)dx\\
&-\beta (n-2)\int_{\mathbb{B}^n}d^{2\beta+\alpha-1}p^{n-3}(\nabla
d\cdot \nabla p) dx.
\end{aligned}
\end{equation}
A direct computation shows that
\[\Delta d=p^2r+\frac{n-1}{r}p\] and

\[\nabla d\cdot \nabla p=p^3r.\] Substituting these above
\begin{equation}
\begin{aligned}
&-\frac{\beta}{\alpha+2\beta}
\int_{\mathbb{B}^n}\text{div}\big(p^{n-2}\nabla
(d^{2\beta+\alpha})\big)dx \\ =& -\beta (2\beta
+\alpha-1)\int_{\mathbb{B}^n}d^{2\beta+\alpha-2}p^n\psi^2dx-(2\beta
+\alpha)\int_{\mathbb{B}^n}
d^{2\beta+\alpha-1}p^n\big(\frac{(n-1)(pr^2+1)}{pr}\big)\psi^2 dx.
\end{aligned}
\end{equation}
We can easily show that \[\frac{pr^2+1}{pr}\ge \frac{1}{d}.\] If
$2\beta+\alpha<0$ then we have

\begin{equation}
-\frac{\beta}{\alpha+2\beta}
\int_{\mathbb{B}^n}\text{div}\big(p^{n-2}\nabla
(d^{2\beta+\alpha})\big)dx \ge
-\beta(2\beta+\alpha+n-2)\int_{\mathbb{B}^n}
d^{2\beta+\alpha-2}p^n\psi^2dx.
\end{equation}
Now we substitute (3.7) into (3.4) and we get

\[\int_{
\mathbb{B}^n}d^{\alpha}|\nabla \phi|^2p^{n-2}\ge
(-\beta^2-\beta(\alpha+n-2)) \int_{
\mathbb{B}^n}d^{2\beta+\alpha-2}\psi^2 p^n dx+\int_{\mathbb{B}^n}
d^{\alpha+2\beta}|\nabla\psi|^2 p^{n-2}dx.\]

Note that the function $\beta\longrightarrow
-\beta^2-\beta(\alpha+n-2)$ attains the maximum for
$\beta=\frac{2-\alpha-n}{2}$, and this maximum is equal to
$(\frac{n+\alpha-2}{2})^2$. Therefore we have the following
inequality

\[\begin{aligned}\int_{
\mathbb{B}^n}d^{\alpha}|\nabla \phi|^2p^{n-2}dx&\ge
\Big(\frac{n+\alpha-2}{2}\Big)^2\int_
{\mathbb{B}^n}d^{\alpha}\frac{\phi^2}{d^2} p^n
dx+\int_{\mathbb{B}^n} d^{2-n}|\nabla\psi|^2 p^{n-2}dx\\
& \ge\Big(\frac{n+\alpha-2}{2}\Big)^2\int_
{\mathbb{B}^n}d^{\alpha}\frac{\phi^2}{d^2} p^n
dx+\int_{\mathbb{B}^n} r^{2-n}|\nabla\psi|^2 dx.
\end{aligned}\]
Notice that the weight function $r^{2-n}$ is in the Muckenhoupt
$A_2$ class and we have the weighted Poincare inequality
\cite{Kenig}( One can also use the reduction of the dimension
technique as in \cite{Brezis-Vazquez}).  Therefore

\begin{equation}
\int_{\mathbb{B}^n}d^{\alpha}|\nabla\phi|^2dx
\ge\Big(\frac{n+\alpha-2}{2}\Big)^2\int_{\mathbb{B}^n}
d^{\alpha}\frac{\phi^2}{d^2}p^ndx+ c\int_{ \mathbb{B}^n}
r^{2-n}\psi^2dx
\end{equation} where $c>0$. Since  $2r\le d\le pr$ then we obtain the following improved Hardy inequality

\begin{equation}\int_{\mathbb{B}^n}d^{\alpha}|\nabla\phi|^2p^{n-2}dx
\ge \Big(\frac{n+\alpha-2}{2}\Big)^2\int_{\mathbb{B}^n}
d^{\alpha}\frac{\phi^2}{d^2}p^ndx+c2^{n-2}\int_{\mathbb{B}^n
}d^{\alpha}\phi^2 dx.\end{equation}

\medskip

It only remains to show that the constant
$(\frac{n+\alpha-2}{2})^2$ is the best constant for the Hardy
inequality (3.8), that is
\[\big(\frac{n+\alpha-2}{2}\big)^2=\inf\Big\{\frac{\int_{\mathbb{B}^n}d^{\alpha}|\nabla \phi|^2 p^{n-2}dx}{\int_{\mathbb{B}^n}d^{\alpha-2}\phi^2p^ndx},
\phi\in C_0^{1}(\mathbb{R}^n), \phi\neq 0\Big\}.\]

Let $\phi_{\epsilon}(d)$ be the family of functions defined by
\begin{equation}
\phi_{\epsilon}(d)=
\begin{cases}
 1 &\quad\text{if}  \quad d \in [0,1],\\
d^{-(\frac{n+\alpha-2}{2}+\epsilon)} &\quad \text{if} \quad d
>1,
\end{cases}
\end{equation}
where $\epsilon>0$ and $d=\text{log}(\frac{1+|x|}{1-|x|})$. It
follows that
\[\int_{
\mathbb{B}^n} d^{\alpha}|\nabla
\phi_{\epsilon}|^2p^{n-2}=\Big(\frac{n+\alpha-2}{2}+\epsilon\Big)^2
\int_{\mathbb{B}^n} d^{-n-2\epsilon} p^ndx\] In the sequel we
indicate $B_1=\{x:d\le 1\}$ $d$-ball centered at the origin in
$\mathbb{B}^n$ with radius $1$.

 By direct computation we get \begin{equation}\begin{aligned}
\int_{\mathbb{B}^n} d^{\alpha}\frac{\phi_{\epsilon}^2}{d^2} p^ndx
&=\int_{B_1}d^{\alpha-2}p^ndx +\int_{\mathbb{B}^n\setminus
B_1}d^{-n-2\epsilon} p^ndx\\
&=\int_{B_1}d^{\alpha-2}p^ndx +
(\frac{n+\alpha-2}{2}+\epsilon)^{-2}\int_{\mathbb{B}^n}d^{\alpha}
|\nabla\phi_{\epsilon}|^2p^ndx.\end{aligned}\end{equation} Since
$n+\alpha-2>0$ then the first integral on the right hand side of
(3.11) is integrable and we conclude by $\epsilon\longrightarrow
0$.

\endproof
\medskip

We now give a new  improved  version of uncertainty principle
inequality on hyperbolic space  which is an immediate consequence
of the improved Hardy inequality (3.9)  and the Cauchy-Schwarz
inequality. Let us mention that a different version of uncertainty
principle inequality on hyperbolic space has been obtained by Sun
\cite{Sun}.

\begin{corollary}(Improved Uncertainty inequality). Let
$\phi\in C_0^{\infty}(\mathbb{B}^n\setminus \{0\})$,
$d=\log(\frac{1+|x|}{1-|x|})$  and $n\ge 2$. Then
\begin{equation}
\Big(\int_{\mathbb{B}^n} d^2 \phi^2 p^n
dx\Big)\Big(\int_{\mathbb{B}^n}|\nabla
\phi|^2p^{n-2}dx-c2^{n-2}\int_{\mathbb{B}^n}\phi^2dx\Big)\ge
\Big(\frac{n-2}{2}\Big)^2\Big(\int_{\mathbb{B}^n}\phi^2 p^n
dx\Big)^2
\end{equation}
where $c>0$.
\end{corollary}
\medskip

 \section*{Improved Weighted Rellich-type inequality}

 Using the same argument as in the proof of Theorem 2.4, we prove
 the following improved Rellich inequality with an explicit
 constant.
\begin{theorem} Let
$\phi\in C_0^{\infty}(\mathbb{B}^n\setminus \{0\})$,
$d=\log(\frac{1+|x|}{1-|x|})$,  $n\ge3$,  $\alpha<2$ and
$n+\alpha-4>0$. Then the following inequality is valid
\begin{equation}
\int_{\mathbb{B}^n}d^{\alpha}|\Delta_{\mathbb{H}} \phi|^2dV \ge
\frac{(n+\alpha-4)^2(n-\alpha)^2}{16}
\int_{\mathbb{B}^n}d^{\alpha}\frac{\phi^2}{d^4}dV+c\epsilon 2^n
\int_{\mathbb{B}^n} d^{\alpha}\phi^2dx
\end{equation}
where $ \epsilon=\frac{(n+\alpha-4)(n-\alpha)}{8}$ and $c>0$.
\proof

A straight forward computation shows that
\begin{equation}
\begin{aligned}
\Delta_{\mathbb{H}}(d^{\alpha-2})&=p^{-n}\text{div}\big(p^{n-2}\nabla(d^{\alpha-2})\big)\\
&=(\alpha-2)(\alpha-3)d^{\alpha-4}+(\alpha-2)(n-1)d^{\alpha-3}(\frac{pr^2+1}{pr}).
\end{aligned}
\end{equation}
Since \[\frac{pr^2+1}{pr}\ge \frac{1}{d}\quad\text {and} \quad
{\alpha}<2,
\] we obtain

\begin{equation}
\Delta_{\mathbb{H}}(d^{\alpha-2})=p^{-n}\text{div}\big(p^{n-2}\nabla(d^{\alpha-2})\big)\le
(\alpha-2)(n+\alpha-4)d^{\alpha-4}.
\end{equation} Multiplying both sides of (3.15) by $\phi^2$ and
integrating,  we obtain

\[\begin{aligned}\int_{\mathbb{H}^n}\Delta_{\mathbb{H}}(d^{\alpha-2})\phi^2dV&=\int_{\mathbb{H}^n}d^{\alpha-2}\Delta_{\mathbb{H}}(\phi^2)dV \\
&=2\int_{\mathbb{H}^n} (\phi\Delta_{ \mathbb{H}}\phi)
d^{\alpha-2}dV+2\int_{\mathbb{H}}|\nabla_{\mathbb{H}^n}\phi|^2d^{\alpha-2}dV\\
&\le(\alpha-2)(n+\alpha-4)\int_{\mathbb{H}^n}d^{\alpha-4}\phi^2dV.
\end{aligned}\]
Therefore

 \begin{equation}
-2\int_{\mathbb{B}^n} (\phi\Delta_{ \mathbb{H}}\phi)
d^{\alpha-2}dV \ge
2\int_{\mathbb{B}^n}|\nabla_{\mathbb{H}}\phi|^2d^{\alpha-2}dV-(\alpha-2)(n+\alpha-4)\int_{\mathbb{B}^n}d^{\alpha-4}\phi^2dV.\end{equation}

Let us apply Young's inequality to expression
$\int_M(\phi\Delta\phi)\rho^{\alpha-2}$ in (3.16) and we obtain
\begin{equation}-\int _{\mathbb{B}^n}\rho^{\alpha-2}\phi\Delta_{\mathbb{H}}\phi dV\le \epsilon\int
_{ \mathbb{B}^n}\rho^{\alpha-4}\phi^2dV+\frac{1}{4\epsilon} \int
_{\mathbb{B}^n}\rho^{\alpha}
|\Delta_{\mathbb{H}}\phi|^2dV\end{equation} where $\epsilon>0$.
Substituting (3.17) into (3.16) and using the improved Hardy
inequality (3.9), we get
\[\frac{1}{4\epsilon}\int _{ \mathbb{B}^n}\rho^{\alpha}
|\Delta_{\mathbb{H}}\phi|^2dV\ge
\Big(\frac{(n+\alpha-4)(n-\alpha)}{4}-\epsilon\Big)\int _{
\mathbb{B}^n}\rho^{\alpha-4}\phi^2dV+c2^{n-2}\int_{\mathbb{B}^n}
d^{\alpha}\phi^2dx\] Since $n+\alpha-4>0$ and $n-\alpha>0$ then we
choose $\epsilon=\frac{(n+\alpha-4)(n-\alpha)}{8}$. Therefore we
obtain the following improved Rellich inequality

\[\int_{\mathbb{H}^n}d^{\alpha}|\Delta_{\mathbb{H}} \phi|^2dV \ge
\frac{(n+\alpha-4)^2(n-\alpha)^2}{16}
\int_{\mathbb{B}^n}d^{\alpha}\frac{\phi^2}{d^4}dV+c\epsilon 2^n
\int_{\mathbb{B}^n} d^{\alpha}\phi^2dx.\] \qed
\end{theorem}

\bibliographystyle{amsalpha}

\begin{thebibliography}{A}

\bibitem [ACP] {Peral}  B. Abdellaoui, D. Colorado, I. Peral, \textit{ Some improved Caffarelli-Kohn-Nirenberg inequalities},
 Calc. Var. Partial Differential Equations {\bf 23}  (2005), no. 3,
 327-345.

\bibitem[BG]{Baras-Goldstein} P. Baras and J. A. Goldstein, \textit
{The heat equation with a singular potential}, Trans. Amer. Math.
Soc. {\bf284} (1984), 121--139

\bibitem [B] {Barbatis} G. Barbatis, \textit { Best constants for higher-order Rellich inequalities in
$L^p(\Omega)$}, preprint

\bibitem [BV]{Brezis-Vazquez} H. Brezis and J. L. V\'azquez, \textit{Blow-up solutions of some nonlinear elliptic
problems}, Rev. Mat. Univ. Complutense Madrid \textbf{10} (1997),
443-469.

\bibitem [CM]{Cabre-Martel} X. Cabr\'e and Y. Martel, \textit{Existence versus explosion instantanée pour des équations de la chaleur linéaires avec potentiel singulier}
, C. R. Acad. Sci. Paris Sér. I. Math., \textbf{329} (1999),
973-978.

\bibitem [C]{Carron} G. Carron, Inégalités de Hardy sur les variétés riemanniennes
non-compactes, J. Math. Pures Appl. (9) 76 (1997), 883-891.

\bibitem [DH] {Davies-Hinz} E. B. Davies, and A. M. Hinz, \textit{Explicit constants for Rellich inequalities in $L\sb p(\Omega)$},
 Math. Z. \textbf {227} (1998), no. 3, 511-523.

\bibitem [D] {Dupaigne} L. Dupaigne, \textit { A nonlinear elliptic PDE with the inverse-square potential
}, J. Anal. Math. \textbf{86}  (2002), 359-398.

\bibitem [FS]{Folland-Sitaram} G. B. Folland and A. Sitaram, \textit{The Uncertainty Principle: A Mathematical Survey},
J. Fourier Anal. Appl. \textbf{3} (1997), 207-238.

\bibitem [FKS] {Kenig} E. Fabes, C. Kenig and R. Serapioni, The local regularity of
solutions of degenerate elliptic equations, Comm. in P.D.E.,
\textbf{7} (1982), 77-116.

\bibitem [GP]{Garcia-Peral}
J. Garcia Azorero and I. Peral, \textit{Hardy inequalities and
some critical elliptic and parabolic problems}, J. Diff.
Equations, {\bf 144} (1998), 441-476.

\bibitem [G] {Grillo} G. Grillo, \textit{ Hardy and Rellich-Type Inequalities for
metrics Defined by Vector fields}, Potential Analysis 18 (2003),
187-217.

\bibitem [GZ] {Goldstein-Zhang} J. A. Goldstein and Qi. S. Zhang, \textit{ Linear parabolic equations with strong singular
potentials}, Trans. Amer. Math. Soc. 355 (2003), 197-211

\bibitem [H]{Heisenberg} W. Heisenberg,  \textit{\"{U}ber den anschaulichen Inhalt der quantentheoretischen Kinematik und Mechanik},
Z. Physik \textit{43} (1927), 172–198.

\bibitem [K1]{Kombe1} I. Kombe,
\textit { The linear heat equation with a highly singular,
oscillating potential}, Proc. Amer. Math. Soc. {\bf 132} (2004),
2683-2691.

\bibitem [K2]{Kombe2} I. Kombe, \textit{ Hardy, Rellich and Uncertainty principle inequalities on Carnot
Groups}

\bibitem [LW]{Li-Wang}P. Li and J. Wang \textit{ Weighted Poincaré inequality and rigidity of complete
manifolds}, preprint

\bibitem [L] {Lindqvist} P. Lindqvist, \textit {On the equation}
$\text{div}(|\nabla u|^{p-2}\nabla u)+\lambda|u|^{p-2}u=0$, Proc.
Amer. Math. Soc. \textbf(109) (1990), 157-164.

\bibitem [PV]{Peral-Vazquez} I. Peral and J. L. V\'azquez,
\textit{ On the stability or instability of the singular solution
of the semilinear heat equation with exponential reaction term},
Arch. Rational Mech. Anal. {\bf 129} (1995),  201-224.

\bibitem [R]{Rellich} F. Rellich, ``Perturbation theory of eigenvalue problems", Gordon and
Breach, New York, 1969.

\bibitem [S]{Sun}  L. Sun, \textit{ An Uncertainty Principle on Hyperbolic
Space}, Proc. Amer. Math. Soc. \textbf{121} (1994), 471-479.

\bibitem [VZ]{Vazquez-Zuazua} J. L. V\'azquez and E. Zuazua, \textit{The Hardy constant and the asymptotic behaviour of the heat equation with an inverse-square potential}
, J. Funct. Anal. \textbf{173} (2000), 103-153.

\end{thebibliography}

\end{document}